\documentclass[12pt]{amsart}
\usepackage{amssymb}

\theoremstyle{plain}
\newtheorem{theorem}{Theorem}

\newtheorem{lemma}[theorem]{Lemma}
\newtheorem{proposition}[theorem]{Proposition}

\theoremstyle{definition}

\theoremstyle{remark}

\def\M{\mathfrak{M}}
\def\Bl1{\mathcal B(\ell_1)}
\def\lin{{\rm lin}}
\def\mid{\::\:}

\newcommand{\abs}[1]{\lvert#1\rvert}
\newcommand{\norm}[1]{\lVert#1\rVert}
\renewcommand{\le}{\leqslant}
\renewcommand{\ge}{\geqslant}

\begin{document} 

\baselineskip=16pt

\title[Lomonosov's theorem]
  {Lomonosov's theorem cannot be extended\\
   to chains of four operators}

\author[V.~G.~Troitsky]{Vladimir~G.~Troitsky}
\address{
         Mathematics Department\\
         University of Illinois at Urbana-Champaign\\
         1409 West Green St.\\
         Urbana, Il 61801\\
         USA}
\email{vladimir@math.uiuc.edu}

\thanks{Supported in part by NSF Grant DMS 96-22454.}

\date{February 12, 1998}

\begin{abstract}
  We show that the celebrated Lomonosov theorem cannot be improved by
  increasing the number of commuting operators. Specifically, we prove
  that if $T\colon\ell_1\to \ell_1$ is the operator without a non-trivial
  closed invariant subspace constructed by C.~J.~Read, then there are
  three operators $S_1$, $S_2$ and $K$ (non-multi\-ples of the
  identity) such that $T$ commutes with $S_1$, $S_1$ commutes with
  $S_2$, $S_2$ commutes with $K$, and $K$ is compact. It is also shown
  that the commutant of $T$ contains only series of $T$.
\end{abstract}

\maketitle

All Banach spaces in this note are assumed to be infinite dimensional
and separable; all operators are linear and bounded.  We say that an
operator is a non-scalar operator if it is not a multiple of the
identity operator.

One of the major results in the history of the Invariant Subspace
Problem was obtained by V.~Lomonosov in~\cite{L} who proved that if an
operator $T$ on a Banach space commutes with another non-scalar
operator $S$ and $S$ commutes with a non-zero compact operator $K$,
then $T$ has an invariant subspace.  Motivated by their study of the
Invariant Subspace Problem for positive operators on Banach lattices,
Y.~A.~Abramovich and C.~D.~Aliprantis have asked recently whether or
not Lomono\-sov's theorem can be extended to chains of four or more
operators.  The purpose of this note is to answer this question in the
negative. For our initial operator $T$ we will take an operator
without an invariant subspace on $\ell_1$ coming from the famous
construction of C.~J.~Read (see~\cite{R1}). Then we will produce three
non-scalar operators $S_1$, $S_2$, and $K$ with $K$ compact (as a
matter of fact $K$ has rank one) such that $TS_1=S_1T$,
$S_1S_2=S_2S_1$, and $S_2K=KS_2$. This will be done in
Section~\ref{s:chain}.

After that, in Section~\ref{s:comm}, we will consider  a related question
of describing the commutant of the C.~J.~Read operator.

\section{A chain from C.~J.~Read's operator to a rank-one operator}
\label{s:chain}

We begin with reminding the reader of the construction in~\cite{R1}
that will be central for us. As in~\cite{R1}, we denote the standard
unit vectors of $\ell_1$ by $(f_i)_{i=0}^\infty$. The symbol $F$
denotes the linear subspace of $\ell_1$, spanned by $f_i$'s, and thus,
$F$ consists of eventually vanishing sequences.

Let ${\bf d}=(a_1,b_1,a_2,b_2,\ldots)$ be a strictly increasing
sequence of positive integers. Also let $a_0=1$, $v_0=0$,
and $v_n=n(a_n+b_n)$ for $n\ge 1$. Then there is a unique sequence
$(e_i)_{i=0}^\infty\subset F$ with the following properties:
\begin{enumerate}
  \item[0)] $f_0=e_0$; 
  \item[A)] if integers $r$, $n$, and $i$ satisfy $0<r\le n$,
    $i\in[0,v_{n-r}]+ra_n$, we have $f_i=a_{n-r}(e_i-e_{i-ra_n})$;
  \item[B)] if integers $r$, $n$, and $i$ satisfy $1\le r<n$,
    $i\in(ra_n+v_{n-r},(r+1)a_n)$, (respectively, $1\le n$,
    $i\in(v_{n-1},a_n))$, then $f_i=2^{(h-i)/\sqrt{a_n}}e_i$, where
    $h=(r+\frac{1}{2})a_n$ (respectively, $h=\frac{1}{2}a_n$);
  \item[C)] if integers $r$, $n$, and $i$ satisfy $1\le r\le n$,
    $i\in[r(a_n+b_n),na_n+rb_n]$, then $f_i=e_i-b_ne_{i-b_n}$;
  \item[D)] if integers $r$, $n$, and $i$ satisfy $0\le r<n$,
  $i\in(na_n+rb_n,(r+1)(a_n+b_n))$, then $f_i=2^{(h-i)/\sqrt{b_n}}e_i$,
  where $h=(r+\frac{1}{2})b_n$.
\end{enumerate}

Indeed, since $f_i=\sum_{j=0}^{i}\lambda_{ij}e_j$ for each $i\ge 0$
and $\lambda_{ii}$ is always nonzero, this linear relation is
invertible. Further,
\begin{equation} \label{e:lins}
  \tag{$*$}
  \lin\{e_i\mid i=1,\dots,n\}=
  \lin\{f_i\mid i=1,\dots,n\}\mbox{ for every }n\ge 0.
\end{equation}
In particular all $e_i$ are linearly independent and also span $F$.
Then C.~J.~Read defines $T\colon F\to F$ to be the unique linear map such that
$Te_i=e_{i+1}$, and in Lemma~5.1 he proves that $\norm{Tf_i}\le1$ for
every $i\ge 0$ provided {\bf d} increases sufficiently rapidly, i.~e.,
satisfies several conditions of the form
\begin{eqnarray*}
  a_n & \ge & G(n, a_0, b_0, a_1, b_1, \dots, a_{n-1}, b_{n-1}),
  \mbox{ and} \\
  b_n & \ge & H(n, a_0, b_0, a_1, b_1, \dots, a_{n-1}, b_{n-1}, a_n),
\end{eqnarray*}
where the $G$ and $H$ are some real-valued functions.
It follows that $T$ can be extended to a bounded operator on $\ell_1$.
Finally, C.~.J.~Read proves that this extension, which is still denoted by
$T$, has no invariant subspaces provided {\bf d} increases
sufficiently rapidly.

Throughout this section we will assume, without loss of generality, that
all integers $a_i$ and $b_i$ are even. We are going to construct non-scalar
operators $S_1$, $S_2$, and $K$ such that
$K$ has rank one and commutes with
$S_2$, $S_2$ commutes with $S_1$, and $S_1$ commutes with $T$. In fact,
we take $S_1=T^2$, so that the equality $TS_1=S_1T$
is automatic. Define $S_2$ on $F$
via
\[
  S_2e_i=
    \left\{
        \begin{array}{ll}
          e_i       & \mbox{ if $i$ is even;}\\
          0         & \mbox{ otherwise. }
        \end{array}
    \right.
\]
We claim that
\[
  S_2f_i=
    \left\{
      \begin{array}{ll}
        f_i       & \mbox{ if $i$ is even;}\\
        0         & \mbox{ otherwise. }
      \end{array}
    \right.
\]
To prove this we consider all possible cases:
\begin{enumerate}
  \item[0)] In this case $S_2f_0=S_2e_0=e_0=f_0$;
  \item[A)] Since $a_n$ is even then
    \[
      S_2f_i=a_{n-r}(S_2e_i-S_2e_{i-ra_n})=
       \left\{
           \begin{array}{ll}
             a_{n-r}(e_i-e_{i-ra_n})=f_i  & \mbox{ if $i$ is even;}\\
             0                            & \mbox{ otherwise. }
           \end{array}
       \right.
    \]
  \item[B)] In this case
    \[
      S_2f_i=2^{(h-i)/\sqrt{a_n}}S_2e_i=
       \left\{
           \begin{array}{ll}
             2^{(h-i)/\sqrt{a_n}}e_i=f_i  & \mbox{ if $i$ is even;}\\
             0                            & \mbox{ otherwise. }
           \end{array}
       \right.
   \]
   \item[C)] Since $b_n$ is even, we have
    \[
      S_2f_i= S_2e_i- b_nS_2e_{i-b_n}=
       \left\{
           \begin{array}{ll}
             e_i-b_ne_{i-b_n}=f_i  & \mbox{ if $i$ is even;}\\
             0                     & \mbox{ otherwise. }
           \end{array}
       \right.
    \]
   \item[D)] Finally, in this case
    \[
      S_2f_i=2^{(h-i)/\sqrt{b_n}}S_2e_i=
       \left\{
           \begin{array}{ll}
             2^{(h-i)/\sqrt{b_n}}e_i=f_i  & \mbox{ if $i$ is even;}\\
             0                            & \mbox{ otherwise.}
           \end{array}
       \right.
   \]
\end{enumerate}

In particular, $S_2$ is bounded on $F$ and can be extended to $\ell_1$.
For every $i\ge 0$ we have
\[
  T^2S_2e_i=
       \left\{
           \begin{array}{ll}
             T^2e_i=e_{i+2} & \mbox{ if $i$ is even;}\\
             0              & \mbox{ otherwise.}
           \end{array}
       \right.
\]
On the other hand
\[
  S_2T^2e_i=S_2e_{i+2}=
       \left\{
           \begin{array}{ll}
             e_{i+2} & \mbox{ if $i$ is even;}\\
             0               & \mbox{ otherwise,}
           \end{array}
       \right.
\]
so that $T^2S_2x=S_2T^2x$ for every $x\in F$. Since $F$ is dense in
$\ell_1$, it follows that $T^2$ and $S_2$ commute on $\ell_1$.

Finally, define $K$ on $\ell_1$ via $Kf_0=f_0$ and $Kf_i=0$ for all
$i>0$.  Then $K$ is a bounded rank one operator on $\ell_1$, and $K$
commutes with $S_2$.

Note that if $m$ divides $a_n$ and $b_n$ for every $n$, then,
similarly to the previous construction, we could take for $S_1$ the
operator $T^m$ instead of $T^2$.  It follows from Lomonosov's theorem
that $T^m$ has an invariant subspace (confer~\cite[Lemma~6.4]{R1}).

In~\cite{R2} C.~J.~Read presents as a modification of his original
example a quasinilpotent operator on $\ell_1$ without closed
nontrivial invariant subspaces. The same argument as above provides a
chain of four commuting operators connecting this operator to a compact
operator.

\section{Commutants of C.~J.~Read's operators}
\label{s:comm}

Let $\M$ denote the collection of all real infinite matrices
$A=(a_{ij})_{i,j=0}^\infty$. Elements of $\M$ can be added entry-wise,
the zero and identity matrices are defined in the natural way. Let
$\mathcal F$ denote the subfamily of $\M$ consisting of the matrices with
finite number of nonzero entries in every column and every row. For
$A$, $B$, and $C$ in $\M$ we say that $AB=C$ if for all $i,j\in\mathbb
N$ we have $a_{ij}=\sum_{k=0}^\infty b_{ik}c_{kj}$ and the series
converge absolutely. Though $AB$ may not exist in general, $AB$ exists
if $A$ or $B$ belongs to $\mathcal F$.

We can also define the action of a matrix on a sequence: we say that
$Ax=y$ for $A\in\M$ and $x,y\in\mathbb R^\mathbb N$ if
$y_i=\sum_{j=0}^\infty a_{ij}x_j$ for every $i\ge 0$. Again, if
$A\in\mathcal F$, then $Ax$ is defined for every $x\in\mathbb
R^\mathbb N$. Let $A^{(j)}$ denote the $j$-th column of $A\in\M$,
then $A^{(j)}=Af_j$. Finally, $(A)_{ij}$ or $a_{ij}$ will denote
the $(i,j)$-th entry of $A\in\M$.

The space $\Bl1$ of all (bounded)  operators on $\ell_1$ can be
naturally embedded in $\M$: if $R\in\Bl1$, then
$r_{ij}=(Rf_j)_i$. Obviously, the sum of operators in $\Bl1$ corresponds
to the sum of matrices in $\M$. Moreover, the action of $R$ on an
element of $\ell_1$ is in accord with the definition of the action of
a matrix on a sequence: if $x=(x_1,x_2,\dots)\in\ell_1$ then
$x=\sum_{j=0}^{\infty}x_jf_j$, so that $Rx=\sum_{j=0}^{\infty}x_jRf_j$,
and
$(Rx)_i=\sum_{j=0}^{\infty}x_j(Rf_j)_i=\sum_{j=0}^{\infty}r_{ij}x_j$.
Also, the product of two operators in $\Bl1$ corresponds to the
product of two matrices: if $R,P\in\Bl1$, then
$RPf_j=RP^{(j)}=R\sum_{k=0}^{\infty}p_{kj}f_k=\sum_{k=0}^{\infty}p_{kj}Rf_k$,
so that $(RP)_{ij}=(RPf_j)_i=\sum_{k=0}^{\infty}p_{kj}r_{ik}$. Notice
that two operators in $\Bl1$ commute if and only if they commute as
matrices in $\M$. Finally, the identity and zero operators correspond
to the identity and zero matrices respectively.

Let $S$ be the right shift operator, i.e., $Sf_j=f_{j+1}$. Given a
formal power series $p(t)=\sum_{n=0}^\infty p_nt^n$, the matrix
$p(S)=p_0I+p_1S+p_2S^2+\dots$ belongs to $\M$ and is of the form
\begin{equation} \label{e:p(S)}
  \left(
    \begin{array}{ccccc}
        p_0  &     0  &     0  &     0  &  \dots \\
        p_1  &   p_0  &     0  &     0  &  \dots \\
        p_2  &   p_1  &   p_0  &     0  &  \dots \\
        p_3  &   p_2  &   p_1  &   p_0  &  \dots \\
      \vdots & \vdots & \vdots & \vdots & \ddots \\
    \end{array}
  \right).
\end{equation}
On the other hand, every matrix of this form is obviously a series
of $S$.

\begin{lemma} \label{l:S-comm}
  If $A\in\M$ commutes with $S$, then $A=p(S)$ for some series $p$.
\end{lemma}

\begin{proof}
  Since
  \[
    (AS)_{ij}=\sum\limits_{k=0}^\infty a_{ik}s_{kj}=a_{i,j+1}
  \]
  and
  \[
    (SA)_{ij}=\sum\limits_{k=0}^\infty s_{ik}a_{kj}=
    \left\{
        \begin{array}{ll}
          a_{i-1,j} & \mbox{ if } i\ge 1;\\
          0         & \mbox{ if } i=0
        \end{array}
    \right.
  \]
  for every pair $(i,j)$, it follows that $A$ is of the form~(\ref{e:p(S)}).
\end{proof}

Consider $Q\in\M$ such that $Qf_j=Q^{(j)}=e_j$ for every $j\ge 0$. It
follows from~(\ref{e:lins}) that $Q\in\mathcal F$. It also follows
from~(\ref{e:lins}) that $Q$ is invertible, and $Q^{-1}\in\mathcal F$.
Further, we can define ``change of basis'' map $A\in\M\mapsto
\tilde{A}=Q^{-1}AQ$. Since $Q,Q^{-1}\in\mathcal F$, this map is defined
for every $A\in\M$, one-to-one, onto, and $A=Q\tilde{A}Q^{-1}$.
Clearly, $\tilde{A}$ describes the action of $A$ in terms of the $e_i$'s.
Finally, $AB=BA$ if and only if
$\tilde{A}\tilde{B}=\tilde{B}\tilde{A}$ for $A,B\in\M$.

Recall that $T$ denotes the Read operator introduced in the previous
section. Since $Te_j=e_{j+1}$ for every $j\ge 0$, then $TQf_j=Qf_{j+1}$,
so that $Q^{-1}TQf_j=f_{j+1}$ which implies $\tilde{T}=S$. Suppose
$RT=TR$ for some $R\in\Bl1$,
then $\tilde{R}S=S\tilde{R}$, so that $\tilde{R}=p(S)$ for some series
$p(t)=\sum_{n=0}^\infty p_nt^n$ by Lemma~\ref{l:S-comm}. Therefore,
$R=Q\tilde{R}Q^{-1}=Qp(S)Q^{-1}=p(QSQ^{-1})=p(T)$.
Since every bounded operator of the form $\sum_{n=0}^\infty
p_nT^n$ commutes with $T$, we have proved the following proposition:

\begin{proposition}
  The commutant of $T$ is the set of all bounded operators of the form
  $\sum_{n=0}^\infty p_nT^n$.
\end{proposition}

I would like to thank Prof.~Y.~A.~Abramovich for suggesting this line
of investigation to me and for our discussions. I am thankful to
Professors C.~D.~Aliprantis, V.~J.~Lomonosov, and C.~J.~Read for their
interest in this work.

We remark in conclusion that in~\cite{TV} we study the modulus $\abs{T}$ of
the Read operator of~\cite{R2}, and we prove that, unlike $T$, the modulus
$\abs{T}$ does have an invariant subspace.

\end{document}